\documentclass[aps,showpacs,pre,amsmath,amsfonts,amssymb,superscriptaddress,nofootinbib,notitlepage,a4paper]{revtex4}
\usepackage{amsthm}

\usepackage{graphicx}
\usepackage{psfrag}
\usepackage[naturalnames]{hyperref}
\usepackage{vmargin}

\setmarginsrb{3 cm}{1.5 cm}{3 cm}{2,5 cm}
{0cm}{1cm}{0cm}{1,2cm} 

\DeclareMathOperator*{\essinf}{ess\,inf}
\def\=d{\overset{\textnormal d}{=}}

\newtheorem{myth}{Theorem}
\newtheorem{mycor}{Corollary}
\newtheorem{myprop}{Proposition}
\newtheorem{mylemma}{Lemma}

\newcommand{\beq}{\begin{equation}} \newcommand{\eeq}{\end{equation}}
\newcommand{\bes}{\begin{split}} \newcommand{\ees}{\end{split}}

\def\w{\omega}

\def\c1t{\widetilde{C}_1}
\def\Csigma{C_{\varsigma}}
\def\rhot{\widetilde{\rho}}
\def\V{\mathcal{V}}
\def\dd{\textrm{d}}
\def\LL{\mathcal{L}}

\begin{document}

\title{The large connectivity limit of the Anderson model on tree graphs}

\author{Victor Bapst}
\affiliation{LPTENS, Unit\'e Mixte de Recherche (UMR 8549) du CNRS et
  de l'ENS, associ\'ee \`a l'UPMC Univ Paris 06, 24 Rue Lhomond, 75231
  Paris Cedex 05, France.}

\begin{abstract}
We consider the Anderson localization problem on the infinite regular tree. Within the localized phase, we derive a rigorous lower bound on the free energy function recently introduced by Aizenman and Warzel. Using a finite volume regularization, we also derive an upper bound on this free energy function. This yields upper and lower bounds on the critical disorder such that all states at a given energy become localized. These bounds are particularly useful in the large connectivity limit where they match, confirming the early predictions of Abou-Chacra, Anderson and Thouless.
\end{abstract}

\maketitle

\section{Introduction}
The strong effect of disorder on the transport properties of a quantum particle was first brought to light by the seminal work of Anderson \cite{Anderson58}. Since then, the understanding of the spectral and dynamical properties of the associated random Schr\"{o}dinger operators has been an ongoing challenge both for mathematicians (monographs include \cite{carmona,PaFi,Ki11}) and for physicists (see \cite{Fifty,Lag09} for recent reviews).
A nice geometry to study this question is that of an infinite regular tree. More specifically, given such a regular rooted tree graph $\mathcal{T}_K$ of fixed branching number $K \geq 2$ (each vertex has $K+1$ neighbors, except the root which has $K$ neighbors), the Anderson model on this graph corresponds to the operator acting on the Hilbert space $\ell^2(\mathcal{T}_K)$ \beq H_{K,t}(\omega) = -t {T}_K+  V(\omega) \ , \eeq
where $T_K$ is the adjacency matrix of $\mathcal{T}_K$, $t >0$ is the hopping strength, and the ``disorder'' $V$ is a random potential, i.e. a set of random variables $V_x$ indexed by the vertices of $\mathcal{T}_K$.

The issue of localization is related to the nature of the spectral measures $\mu_{K,t,\delta_x}(dE,\omega)$ associated with $H_{K,t}(\omega)$ and the Kronecker function $\delta_x \in \ell^2(\mathcal{T}_K)$ localized at $x \in \mathcal{T}_K$. More precisely, the different spectra of $H_{K,t}(\omega)$ are associated with the components of the Lebesgue decomposition of the measures $\mu_{K,t,\delta_x}(dE,\omega)$: pure point, singular continuous or absolutely continuous. Ergodicity and equivalence of the local measures \cite{JaLA20} ensures that the supports of the latter are almost surely non-random \cite{carmona,PaFi,Ki11} and do not depend on $x \in \mathcal{T}_K$. The operator $H_{K,t}(\omega)$ exhibits \textit{spectral localization} in an interval $I \subset \mathbb{R}$ if the spectral measures $\mu_{K,t,\delta_x}(dE,\omega)$ associated with $x \in \mathcal{T}_K$ are almost surely all of only pure point type in $I$. A stronger statement \cite{Ki11,Stoll} is to say that $H_{K,t}(\omega)$ exhibits \textit{exponential dynamical localization} in $I$: roughly speaking, this means that a state initially localized in space and in energy (within $I$), and evolving according to the Schr\"{o}dinger dynamics generated by $H_{K,t}(\omega)$, has a probability of being found a time $\tau$ later at a distance $R$ that decays exponentially fast with $R$, uniformly in $\tau$. An opposite behavior is absolutely continuous spectrum, which is equivalent to the statement that current injected through a wire at the corresponding energy on a site of the lattice will be conducted through the graph to infinity \cite{MiDe94,ASW06}. In the particular case of the regular tree, this transport happens ballistically \cite{AW12_ballistic}.

The first analysis of Anderson localization on the infinite tree was performed by Abou-Chacra, Anderson and Thouless \cite{ACTA73,ACT74}, later followed by many works from the physics community, see \cite{Efetov83,MiFy91,MiDe94,SaBe03,MoGa09,MoGa11,BST10,BaSe11,JoBA12,Thouless74,BTT} and references therein. Besides the opposition between extended and localized spectrum, recent topics of interest concern the nature of the extended and localized phases and of the delocalization transition. From a mathematical point of view, the main results include analyticity of the density of states for random potentials close to the Cauchy distribution \cite{AcKl92}, spectral and dynamical localization of the spectrum for sufficiently strong disorder \cite{AM93,Ai94}, and the persistence of extended states at sufficiently low disorder \cite{Kl98,ASW06,FrHaSp07}. More recently, a major progress has been made by Aizenman and Warzel \cite{aw_long_v3} by providing a criterion for the presence of extended states, which is (almost) complementary to the criterion for localization previously found in \cite{AM93,Ai94}. This allows them to work out the small disorder behavior to a greater extent, unveiling in particular the (unexpected) absence of localized states at any energy within the spectrum for bounded random potentials and sufficiently small disorder \cite{AW11,aw_prl_11}.

An interesting limit in which the phase diagram of the Anderson model can be computed more explicitly is the large connectivity one $K \rightarrow \infty$. This case has been considered early, see the 1974's paper by Thouless \cite{Thouless74} for a review of the different predictions available at that time. However, since then, no progress had been made and this large $K$ limit was still considered unresolved \cite{BST10,AW11}. The main contribution of this paper is to confirm the result obtained by Abou-Chacra, Anderson and Thouless for particular random potentials in \cite{ACTA73}, but on a rigorous basis and for generic disorder distributions. The exactness of this result is somehow surprising since the peculiar nature of the delocalized states near the transition \cite{aw_long_v3,AW11,aw_prl_11,BTT} was not understood at that time. Note that this asymptotic is also the one that can be obtained using the heuristic ``best estimate'' of Anderson's original paper \cite{Anderson58}. Our proof is essentially based on the criterion for extended states given in \cite{aw_long_v3}. Let us finally mention that this large connectivity limit is not only of academic interest, since it has recently been argued to give a good approximation of the localization transition of many-particles systems (the connectivity being related to the number of elementary excitations, and the regular tree an approximation of the Fock space of the system) \cite{AGKL97,AKR10,LuSa12}.
\\

We now turn to the presentation of the result derived in this article. We will assume that the random variables $V_x$ are independent and identically distributed, with a probability density $\rho$ that satisfies:
\begin{description}
\item[(A)] $\rho$ is Lipschitz continuous, with a Lipschitz constant $C$,
\item[(B)] There exist $\varsigma \in (0,1)$ and $\Csigma >0$ such that for all $v \in \mathbb{R}^* = \mathbb{R} \setminus \{0\}$: \beq \rho(v) \leq \frac{\Csigma}{|v|^{1+\varsigma}} \ . \eeq
Under Assumption \textbf{(A)}, it follows that $\varsigma$ could be taken as close to $0$ as desired; but we shall not use this in the following.
\item[(C)] If $\textrm{supp}(\rho)$ is unbounded, then for all $k < \infty$, $\inf_{|v|\leq k} \rho(v) >0$.
\end{description}
These conditions are for instance satisfied by linear combinations of Gaussian or Cauchy distribution. They do not apply to piecewise constant distributions, even though we believe that this is only a limitation of our proof technique. It can easily be seen that they imply the conditions A to E of \cite{aw_long_v3} (our last condition is stated only for this purpose and we will not use it in the following). Standard ergodicity results \cite{carmona,PaFi,Ki11} also imply that the spectrum of $H_{K,t}(\omega)$ is almost surely given by the set sum of the support of $V$ and of that of $-t T_K: \; [-2 t \sqrt{K}, 2t \sqrt{K}]$.
Our result then reads as follows:
\begin{myth} \label{th_clear}
Let $\rho$ satisfy assumptions \textbf{(A-C)},  $E \in \mathbb{R}$ be such that $\rho(E)>0$, $t_K = \frac{g}{K \log K}$ and $g_c(E) = \frac{1}{4 \rho(E)}$. Then:
\begin{itemize} 
\item if $g > g_c(E)$, the operator $H_{K,t_K}(\omega)$ has, for $K$ large enough, almost surely some continuous spectrum in a neighborhood of the energy $E$;
\item if $g < g_c(E)$, there exists $\delta_E>0$ such that for $K$ large enough the spectrum of $H_{K,t_K}(\omega)$ is \textit{exponentially dynamically localized} in the range $[E-\delta_E,E+\delta_E]$, meaning that there exist $\mu>0$ and $C < \infty$ such that $\sum_{\substack{y \in \mathcal{T}_K \\ \textrm{dist}(x,y) = L}} \mathbb{E} \left( \sup_{\tau \in \mathbb{R}} \left | \langle \delta_x , P(H_{K,t}(\w))e^{i \tau H_{K,t}(\w)} \delta_y \rangle \right|^2\right) \leq C e^{-\mu(I) L}$, where $P$ is the spectral projector on $[E-\delta_E,E+\delta_E]$, and $\mathbb{E}$ represents the average with respect to $\w$.
\end{itemize}
\end{myth}
Our result also bears the following corollary: define $t_c(K)$ as the hopping below which only localized spectrum can be found. Then we have the following bounds on $t_c(K)$:
\begin{mycor} \label{cor} Let $\rho$ satisfy assumptions \textbf{(A-C)}. Then, for all $\epsilon >0$ and $K$ large enough: \beq \label{eq_upper_bound_cor}  \frac{1-\epsilon}{K \log K} \frac{1}{4 \| \rho \|_\infty}  \leq t_c(K) \leq \frac{1+\epsilon}{K \log K} \frac{1}{4 \| \rho \|_\infty} \ . \eeq
\end{mycor}
The lower bound improves on the previously best known result which, as explained in \cite{warzel12}, follows from the analysis of \cite{AM93,Ai94}: for all $\epsilon >0$ and $K$ large enough, $t_c(K) \geq \frac{1-\epsilon}{K \log K} \frac{1}{(2e) \| \rho \|_\infty}$. Our result thus improves on this bound by a factor $2/e$, and shows that this lower bound is tight by deriving a (previously unknown) upper bound on $t_c(K)$. Moreover, it also gives a ``local'' result (as stated in Theorem~\ref{th_clear}) which was out of reach with the previous methods.

Note that in this regime where the hopping term $t_K$ is asymptotically much smaller than $1/\sqrt{K}$, the density of states of the model converges in the $K \rightarrow \infty$ limit towards the density of disorder $\rho$. Also note that, as observed in \cite{aw_long_v3}, this result also applies to the corresponding operator on the fully regular tree graph (or Bethe lattice) in which \textit{every} vertex has $K+1$ neighbors.
\\

The plan of our paper is as follows: in the next section we first state and discuss in more details the result that we obtain before writing down an explicit expression for the ``free energy'' function introduced in \cite{aw_long_v3}. In Section~\ref{sec_bounds} we derive upper and lower bounds on this quantity that match when the connectivity $K+1$ becomes large, allowing to prove our result. Section~\ref{sec_numerical} presents a numerical confirmation of our result, while we draw our conclusions in Sec.~\ref{sec_ccl}. Some more technical computations are deferred to an appendix.

\section{The ``free energy'' function $\varphi_{K,t}$ and its computation}

\subsection{Strategy of the proof}
Following \cite{AM93,aw_prl_11, aw_long_v3}, we will be interested in the asymptotic behavior of the Green function between the root $0$ and a site at distance $L$ that we denote $L$; the latter is defined by:
\beq \label{eq_def_G} G_{K,t}(0,L,E+i\eta,\omega) = \left \langle \delta_0, (H_{K,t}(\omega)-E-i\eta)^{-1} \delta_L \right \rangle \ , \eeq
where $\eta  \in \mathbb{R}_+^*$.
Introducing, for $s<1$ \beq \label{eq_def_varphi} \varphi_{K,t}(s,E) = \lim_{L \rightarrow \infty} \lim_{\eta \searrow 0} \frac{ \log \mathbb{E} | G_{K,t}(0,L,E+i\eta,\omega)|^s}{L} \ , \eeq
where $\mathbb{E}$ denotes the average with respect to $\omega$ and the limits are well defined for almost all $E$, and commute \cite{aw_long_v3}, 
and $\varphi_{K,t}(1,E) = \lim_{s \nearrow 1} \varphi_{K,t}(s,E)$ (the limit is again well defined, see \cite{aw_long_v3}), we can state a weak version of the main results of \cite{AM93,Ai94,aw_long_v3}:
\begin{myth}[\cite{AM93,Ai94,aw_long_v3}] \label{th_aw} Assume that $\rho$ satisfies assumptions \textbf{(A-C)}. Let $I \subset \mathbb{R}$ be a fixed interval. Then:
\begin{itemize}
\item if for Lebesgue almost all $E \in I$, $\varphi_{K,t}(1,E) > - \log K$, the spectrum of $H_{K,t}(\omega)$ is almost surely absolutely continuous in $I$, 
\item if for Lebesgue almost all $E \in I$, $\varphi_{K,t}(1,E) < - \log K$, the spectrum of $H_{K,t}(\omega)$ is almost surely of pure point type in $I$. Moreover if $\mathrm{ess} \sup_{E \in I} \varphi_{K,t}(1,E) < - \log K $, then $H_{K,t}(\omega)$ almost surely exhibits exponential dynamical localization in $I$.
\end{itemize}
\end{myth}
Our first result is as follows:
\begin{mylemma} \label{th_precise_1}
Let $\rho$ satisfy assumptions \textbf{(A-C)}, $E \in \mathbb{R}$ be such that $\rho(E)>0$, $t_K = \frac{g}{K \log K}$, and $g_c(E) = \frac{1}{4 \rho(E)}$. Then for all $K$ at which the operator $H_{K,t_K}(\omega)$ exhibits spectral localization (i.e. almost surely pure point spectrum) in a neighborhood $\V$ of $E$, if $g > g_c(E)$ and $K$ is large enough, there exists $\delta_{E,K}>0$ such that $\varphi_{K,t}(1,E') > - \log K$ for almost all $E' \in [E-\delta_{E,K},E+\delta_{E,K}]$
\end{mylemma}
Using Theorem~\ref{th_aw}, this shows that if $g>g_c(E)$ and $K$ is large enough, assuming spectral localization is not consistent and $H_{K,t_k}(\w)$ must thus have some continuous spectrum in any neighborhood of $E$, leading to the first part of Theorem~\ref{th_clear}. The proof of Lemma~\ref{th_precise_1} will go through the computation of the free energy function $\varphi_{K,t}(s,E)$ assuming pure point spectrum in a neighborhood of $E$; we explain how this is done in the next subsection~\ref{sec_comp_loc}. The derivation of the lower bound on $\varphi_{K,t}(s,E)$ is then presented in Sec.~\ref{sec_lower_bound}.

If we had additional regularity results on the free energy function $\varphi_{K,t}(1,E)$, for instance if there was a proof that it is differentiable in $(t,E)$ with only isolated critical point (as suggested in \cite{aw_long_v3}), we could exclude the possibility that $E \rightarrow \varphi_{K,t}(1,E) - \log K$ vanishes (except on isolated points), and replace ``some continuous spectrum'' by ``only absolutely continuous spectrum'' in Theorem~\ref{th_clear}. 
\\

\begin{figure}[h]
\center
\includegraphics[width=12cm]{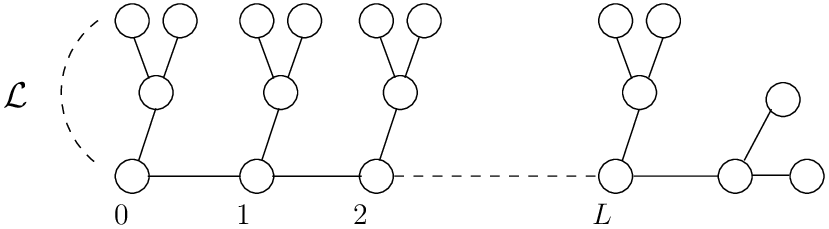}
\caption{A graphical representation of $\Lambda_{L,\LL}$ with $\LL = 2$.}
\label{fig_lambda_L}
\end{figure}

In order to prove localization, it will be more convenient to use a slightly different version of the second part of Theorem~\ref{th_aw} above. We will trade the finite $\eta$ regularization in (\ref{eq_def_G}) for a finite volume one. More precisely, we first denote $1, 2, \dots, L-1$ the sites on the unique path between the root $0$ and site $L$. We then introduce the pre-limit:
\beq \varphi_{K,t}^{(\LL)}(s,E,L) = \frac{ \log \mathbb{E} | G_{K,t}^{(\LL)}(0,L,E,\omega)|^s}{L} \ , \hspace{0.5 cm} G_{K,t}^{(\LL)}(0,L,E,\w) = \left \langle \delta_0, (H_{K,t}^{(L,\LL)}(\omega)-E)^{-1} \delta_L \right \rangle \ , \eeq
where $H_{K,t}^{(L,\LL)}$ denotes the restriction of $H_{K,t}(\w)$ to $\ell^2(\Lambda_{L,\LL})$, with $\Lambda_{L,\LL} = \{ x \in \mathcal{T}_K, \inf_{i \in \{0, \dots, L\}} \textrm{dist}(x,i) \leq \mathcal{L}-1 \}$. A graphical representation of $\Lambda_{L,\LL}$ is shown on Fig.~\ref{fig_lambda_L}. Because the operator $H_{K,t}^{(L,\LL)}(\w)$ acts on an Hilbert space of finite dimension and the parameter $s<1$ regularises the poles of the Green function, the regularizing parameter $\eta$ introduced in (\ref{eq_def_G}) can be directly taken to zero here. Our second result is then as follows:
\begin{mylemma} \label{th_precise_2}
Let $\rho$ satisfy assumptions \textbf{(A-C)}, $E \in \mathbb{R}$ be such that $\rho(E)>0$, $t_K = \frac{g}{K \log K}$, and $g_c(E) = \frac{1}{4 \rho(E)}$. Then  if $g < g_c(E)$ and $K$ is large enough, there exist $s<1$ and $\epsilon_0 >0, \ \delta_{E} >0$ such that for all $\LL \geq 2$ and $L$ large enough, $\varphi_{K,t}^{(\LL)}(s,E',L) < - \log K - \epsilon_0$ for all $E' \in [E-\delta_{E},E+\delta_{E}]$.
\end{mylemma}
The lemma entails the following bound: there exist $s<1, A>0, \delta_E>0, \epsilon_0>0$ such that, for any $L \geq 0$ and $\mathcal{L} \geq 2$
\begin{equation} \int_{E-\delta_E}^{E+\delta_E}  \mathbb{E} \left[ \left \langle \delta_0, (H_{K,t}^{(L,\mathcal{L})}(\omega)-E')^{-1} \delta_L \right \rangle \right]^s   \textrm{d}E' \leq A e^{(-\log K - \epsilon_0) L} \ . \end{equation}
Using Theorem 1.2 of Aizenman's original paper \cite{Ai94}, and the convergence of $H_{K,t}^{(L,\LL)}$ towards $H_{K,t}$ in the strong resolvent sense when $\LL \rightarrow \infty$, this implies that, for $K$ large enough, $H_{K,t}(\omega)$ almost surely exhibits exponential dynamical localization in $[E-\delta_E,E+\delta E]$, and thus implies the second part of Theorem~\ref{th_clear}. We explain how to compute $\varphi_{K,t}^{(\LL)}$ in Sec.~\ref{sec_finite_volume}, and derive an upper bound for it in Sec.~\ref{sec_upper_bound}.
\\
The common technical tool in the proof Lemma~\ref{th_precise_1} and Lemma~\ref{th_precise_2} will be that all the energies and local Green functions involved in the computation of the free energy function will be real, either because one assumes spectral localization, or thanks to the finite volume regularization.

Finally, let us note that the bounds that we shall obtain on $\varphi_{K,t}(1,E)$ and $\varphi_{K,t}^{(\LL)}(s,E,L)$ could in principle be used to deduce lower and upper bounds on a mobility edge (a curve in the $(t,E)$ plane separating localized from extended states) also for finite $K$; however, without further work to optimize them, the latter are very far away from one another unless $K$ becomes extremely large, and therefore we did not pursue this direction.

\subsection{Computation of $\varphi_{K,t}(s,E)$ in the localized phase}
\label{sec_comp_loc}
We first explain the strategy to prove Lemma~\ref{th_precise_1}, and explain how to compute the infinite-volume function $\varphi_{K,t}(s,E)$ in the localized phase, for $s$ close enough to $1$. To do this, let us first explain how one can compute the Green function (\ref{eq_def_G}) between the sites $0$ and $L$. We recall that we denoted $1,2 \dots, L-1$ the sites on the (unique) path between $0$ and $L$; we also denote $L+1$ one of the neighbors of site $L$ different from $L-1$. Given a site $x$, the set of its neighbors $y$ such that $\textrm{dist}(0,y) = \textrm{dist}(0,x)+1$ will be denoted $\mathcal{N}^+_x$. We will use the following expression for $G_{K,t}(0,L,E+ i \eta,\omega)$ (see for example \cite{aw_long_v3}):
\beq \label{eq_green_0_L} G_{K,t}(0,L,E+i\eta,\omega) = \frac{1}{t} \prod_{x=0}^L \frac{t}{V_x(\omega)-E-i\eta- t^2 \sum_{\substack{y\in \mathcal{N}_x^+\\ y \neq x+1}} \Gamma_{K,t}(y,E+i\eta,\omega)- t^2\Gamma_{K,t}(x+1,E+i\eta,\omega)} \ . \eeq 
This expression is standard and can be derived using resolvent identities. In this equation $\Gamma_{K,t}(y,E+i\eta,\omega)$ is the Green function on site $y$ and at energy $E+i\eta$ for the operator $\widetilde{H}_{K,t}$ obtained from $H_{K,t}$ by removing the path between the root $0$ and $y$. It satisfies the recursion equation  \cite{ACTA73,Kl98,BaSe11,aw_long_v3}:
\beq \label{eq_gamma_recursive} \Gamma_{K,t}(y,E+i\eta,\omega) =   \frac{1}{ V_y(\omega) - E-i \eta - t^2 \sum_{z \in \mathcal{N}^+_y} \Gamma_{K,t}(z,E+i\eta, \omega) } \ . \eeq
This is a random quantity because of the randomness in the $V_y$. If we introduce the probability distribution $p_{K,t,E+i\eta}$ of  $\Gamma_{K,t}(y,E+i\eta,\omega)$ with respect to the underlying probability measure, we arrive at the following \textit{recursive distributional equation}:
\beq \label{eq_rde_gamma_im} \Gamma \=d  \frac{1}{ V- E-i \eta - t^2  \sum_{i=1}^K \Gamma_i } \ . \eeq
In this equation, both $\Gamma$ on the left hand side and the $\Gamma_i$ on the right hand side are independent random variables drawn from $p_{K,t,E+i\eta}$, while $V$ is distributed with density $\rho$. This equation is known to admit a unique solution as long as $\eta >0$  \cite{BoLe08}.
\\

In the localized phase (i.e. if the spectrum is almost surely pure point in some neighborhood $\mathcal{V}$ of $E$), for almost all realizations of the disorder, $\lim_{\eta \searrow 0} \Gamma_{K,t}(y,E'+i\eta,\omega)$ exists and is real for all sites $y$ and almost all $E' \in \mathcal{V}$. Conversely, for almost all $E$ within the localized phase, and almost all realizations of the disorder, $\lim_{\eta \searrow 0} \Gamma_{K,t}(y,E+i\eta,\omega)$ exists and is real for all sites $y$ \cite{footnote_1}. Therefore, for an energy $E$ such that the above property holds, the distribution of $\Gamma_{K,t}(y,E+i\eta,\omega)$ becomes supported on the real axis, and (\ref{eq_rde_gamma_im}) can be replaced by:
\beq \label{eq_rde_gamma_real} \Gamma \=d  \frac{1}{ V - E - t^2 \sum_{i=1}^K \Gamma_i } \ , \eeq
where $\Gamma$ is now real, and the solution of this equation has to be selected as the limit of the one of (\ref{eq_rde_gamma_im}) when $\eta \searrow 0$. Since the random variable $V$ admits a density, its sum with any other random variable also does \cite{feller1}, and we deduce from (\ref{eq_rde_gamma_real}) that $\Gamma$  admits a probability density, that we denote $p_{K,t,E}$.

Because $s<1$ the expectation value of  $|G(0,L,E+i \eta,\omega)|^s$ conditioned on the values of the disorder on all the sites $ x \notin \{0, 1, \dots, L\}$ is uniformly bounded with respect to this disorder (this argument can be found in \cite{AM93, aw_long_v3})
\beq \label{eq_bound_power_s} \mathbb{E} \left[ |ÊG_{K,t}(0,L,E+i\eta,\omega)|^s | \{V_x\}_{x \notin \{0, \dots,L\}}\right] \leq \left(2^s \frac{\|\rho\|^s_\infty}{1-s} \right)^{L+1} t^{sL} \ . \eeq This can be obtained here by noting that only the first $x$ factors in the product defining (\ref{eq_green_0_L}) depend on $V_x$. (\ref{eq_bound_power_s}) then follows from averaging the absolute value of (\ref{eq_green_0_L}) raised to the power $s$ over the $V_x$ from $x=0$ up to $x=L$, and using repeatedly the bound $\mathbb{E} \left [ 1/|V-a|^s \right] \leq 2^s \frac{\|\rho\|^s_\infty}{1-s}$ (cf. Eq. (A5) of \cite{aw_long_v3}). It follows that, using the convergence in distribution of $\Gamma_{K,t}(y,E+i\eta)$ towards $\Gamma_{K,t}(y,E)$, we can substitute $\Gamma_{K,t}(y,E)$ for $\Gamma_{K,t}(y,E+i\eta)$ in  $\lim_{\eta \searrow 0} \mathbb{E} | G_{K,t}(0,L,E+i \eta,\omega)|^s$. Introducing the probability density $\rho_{K,t,E}$ of $V-E - t^2 \sum_{i=1}^{K-1} \Gamma_i$ (with the $\Gamma_i$ independently distributed according to the solution $p_{K,t,E}$ of (\ref{eq_rde_gamma_real})), we arrive at the expression that will be our starting point for the following study:
\beq \begin{split} \label{eq_G_unfolded_final} \lim_{\eta \searrow 0} \mathbb{E} | G_{K,t}(0,L,E+i \eta,\omega)|^s &= \frac{1}{t^s} \int  \prod_{x=0}^L \left[ \frac{t^s}{|e_x-t^2 \Gamma_{K,t}(x+1,E)|^s} \rho_{K,t,E}(e_x) \dd e_x \right] \\ & \hspace{3.5 cm}p_{K,t,E}(\Gamma_{K,t}(L+1,E)) \dd \Gamma_{K,t}(L+1,E) \ , \end{split} \eeq
where we recall that, with these notations, 
\beq \Gamma_{K,t}(x,E) = \frac{1}{e_{x}-t^2 \Gamma_{K,t}(x+1,E)} \ .  \eeq
In particular $\Gamma_{K,t}(x,E)$ depends on (and only on) all the ``energies'' $(e_{y})_{x \leq y \leq L}$, and on $\Gamma_{K,t}(L+1,E)$. We now introduce \beq \label{eq_def_F_infinite}  \begin{array}{l l l l}F_{K,t,s,E}: &L^1(\mathbb{R}) & \longrightarrow \;\; & L^1(\mathbb{R})  \\  &a(x)& \longrightarrow &\int_{-\infty}^\infty \mathcal{K}_{K,t,s,E}(x,y) a(y) \dd y \ , \end{array} \eeq 
with \beq \label{eq_def_kernel} \mathcal{K}_{K,t,s,E}(x,y) = \frac{t^{2-s}}{|x|^{2-s}} \rho_{K,t,E}\left(-\frac{t^2}{x}-y\right) \ .\eeq
Then it follows from (\ref{eq_G_unfolded_final}) and changes of variables $e_x \rightarrow \frac{-t^2}{e_x-t^2 \Gamma_{K,t}(x+1,E)}$ that:
\beq \lim_{\eta \searrow 0} \mathbb{E} | G_{K,t}(0,L,E+i \eta,\omega)|^s  = \frac{1}{t^s} \| F_{K,t,s,E}^{L+1} (p_{K,t,E}(- \cdot/t^2)) \|_1  \ , \eeq
where the initial condition is given in terms of the probability density $p_{K,t,E}$ of the solution $\Gamma$  of (\ref{eq_rde_gamma_real}).
Finally, we arrive at: \beq \label{eq_varphi_final} \varphi_{K,t}(1,E)  = \lim_{s \nearrow 1} \lim_{L \rightarrow \infty} \log \left( \| F_{K,t,s,E}^{L} (p_{K,t,E}(-\cdot/t^2)) \|^{\frac{1}{L}}_1 \right) \ , \eeq
which, we recall, is valid for all $E$ such that $\lim_{\eta \searrow 0} \Gamma_{K,t}(y,E+i\eta,\omega)$ exists and is real for all sites $y$ and almost all $\omega$.

\subsection{Computation of the finite-volume function $\varphi_{K,t}^{(\LL)}(s,E,L)$}
\label{sec_finite_volume}
We now explain how to compute the finite-volume free energy $\varphi_{K,t}^{(\LL)}(s,E,L)$. We follow exactly the same steps as before, only substituting $H_{K,t}^{(L,\LL)}$ for $H_{K,t}$ everywhere. In this case, thanks to the finite-volume regularization, we do not have to worry about the $\eta \searrow 0$ limit, and the local resolvents $\Gamma_{K,t}^{(\LL)}$ (defined as before, and whose distribution does \textit{not} depend on $L$) are automatically real. The only difference appears in their probability density, which must now be computed with $H_{K,t}^{(L,\LL)}$ instead of $H_{K,t}$. In particular, the density $p_{K,t,E}^{(\LL)}$ of $\Gamma_{K,t}^{(\LL)}(y,E,\w) \equiv \Gamma^{(\LL)}$ is obtained by $\LL$ iterations of the recursive equation (\ref{eq_rde_gamma_real}):
\beq \Gamma^{(j+1)} \=d  \frac{1}{ V - E - t^2 \sum_{i=1}^K \Gamma_i^{(j)} } \ ,   \eeq
where the $\Gamma_i^{(j)}$ on the right hand side are independently drawn from $p_{K,t,E}^{(j)}$, with the initial condition $\Gamma^{(1)} \=d \frac{1}{V-E}$.  We define the density $\rho_{K,t,E}^{(\LL)}$ of $V-E-t^2 \sum_{i=1}^{K-1} \Gamma_i^{(\LL)}$ accordingly, and obtain $\varphi_{K,t}^{(\LL)}(s,E,L)$ as:
\beq \varphi_{K,t}^{(\LL)}(s,E,L) = \log \left( \| \left(F^{(\LL)}_{K,t,s,E}\right)^L (p_{K,t,E}^{(\LL)}(-\cdot/t^2)) \|^{\frac{1}{L}}_1 \right) \eeq
where we defined, as in Eq. (\ref{eq_def_F_infinite}-\ref{eq_def_kernel}): \beq  \begin{array}{l l l l}F^{(\LL)}_{K,t,s,E}: &L^1(\mathbb{R}) & \longrightarrow \;\; & L^1(\mathbb{R})  \\  &a(x)& \longrightarrow &\int_{-\infty}^\infty \mathcal{K}^{(\LL)}_{K,t,s,E}(x,y) a(y) \dd y \ , \end{array} \eeq 
with: \beq \mathcal{K}^{(\LL)}_{K,t,s,E}(x,y) = \frac{t^{2-s}}{|x|^{2-s}} \rho^{(\LL)}_{K,t,E}\left(-\frac{t^2}{x}-y\right) \ .\eeq

\subsection{Connection with the criterion of Abou-Chacra, Anderson and Thouless and sketch of proof}
\label{sec_connection}
In the next section, we shall compute the sign of $\varphi_{K,\frac{g}{K\log K}}(1,E)+ \log K$ and $\varphi_{K,\frac{g}{K\log K}}^{(\LL)}(s,E,L)+ \log K$ when $g$ is fixed, $K$ is large. But before that, let us briefly explain heuristically how we can recover the critical condition given by Abou-Chacra, Anderson and Thouless  in their paper \cite{ACTA73} (see also \cite{MiFy91} for an alternative derivation) from this point.
\\

Let us proceed as if the operator $F_{K,t,s,E}$ was a finite dimensional one. Then using the positivity of the kernel $\mathcal{K}_{K,t,s,E}$, we could apply the Perron-Frobenius theorem and conclude that $\varphi_{K,t}(s,E) = \log \lambda_{K,t,s,E}$ where $\lambda_{K,t,s,E}$ is the largest eigenvalue of the operator $F_{K,t,s,E}$.  The critical condition would then correspond to $\lambda_{K,t,1,E} =1/K$, and the critical value of the hopping would be given by the largest value of $t$ for which the equation
\beq \label{eq_ACTA} a(x) = \frac{Kt^{2-s}}{|x|^{2-s}} \int \rho_{K,t,E} \left( -y - \frac{t^2}{x} \right) a(y) \dd y \eeq possesses a solution up to $s =1$. Making the substitutions $t \rightarrow V$, $\rho_{K,t,E}(z) \rightarrow Q(-z)$ and $s \rightarrow 2(1-\beta)$, we obtain Eq. (6.7) of \cite{ACTA73}. A closer look at our derivation of (\ref{eq_varphi_final}) reveals that the eigenvector of the largest eigenvalue of (\ref{eq_ACTA}), $a(x)$, has a simple interpretation: it is, up to a normalization, the limiting conditional expectation value of $|G(0,L,E+i\eta, \omega)|^s$, conditioned on $\Gamma_{K,t}(1,E)$ being equal to $-x/t^2$. A very similar interpretation can be given for the dominant eigenvector of the adjoint equation of (\ref{eq_ACTA}), used in \cite{ACTA73}. Note that the interpretation that we obtain for these eigenvectors is somehow more explicit than the one given in \cite{ACTA73}, where they appear as being related to the amplitude of the tails of the distribution of the local Green function $G(0,0,E+i\eta,\omega)$, assuming a power-law tail with an exponent $1+\beta$.

\begin{figure}[h]
\includegraphics[width=6.5cm]{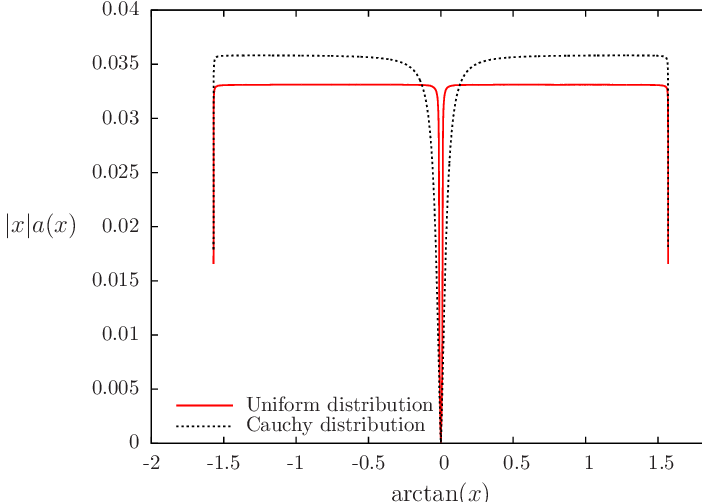} \hspace{0.5 cm}
\includegraphics[width=6.5cm]{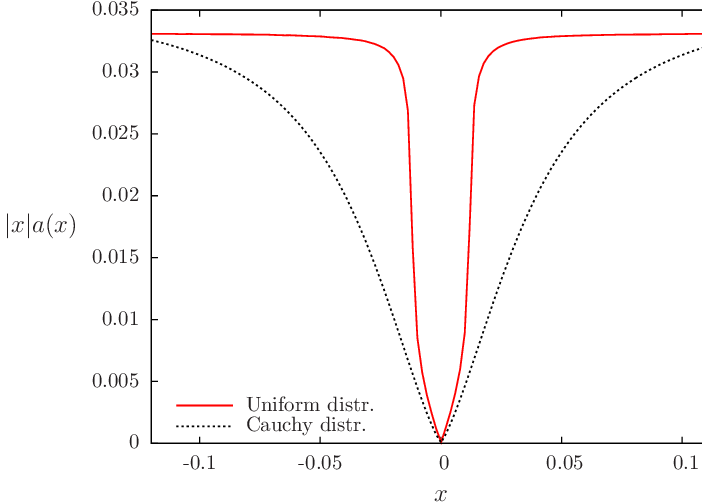}
\caption{$|x|a(x)$, where $a(x)$ is the eigenvector of largest eigenvalue for a discretization of the kernel $\mathcal{K}_{K,t,s,E}$, for $K=2$, $s=1$, $E=0$, and two densities of disorder $\rho$: a uniform distribution $\rho(v) = \frac{1}{2}\mathbf{1}(v \in [-1,1])$ shown with a solid black line (for $t=0.11$), and a standard Cauchy distribution $\rho(v) = \frac{1}{ \pi(v^2+1)}$ shown with a dotted red line (for $t=0.23$). Right: detail around $x=0$. $a(x)$ has a finite limit when $x \rightarrow 0$. In each case $t$ is chosen close to its critical value (the precise value of $t$ has no influence on the global shape of the eigenvector). $a(x)$ was obtained by iterations of the kernel $\mathcal{K}_{K,t,s,E}$, using a discretization of the real line in 16000 elementary intervals.}
\label{fig_eigenvector}
\end{figure}

Here it is worth showing the numerically determined eigenvector of largest eigenvalue of a finite dimensional approximation of $\mathcal{K}_{K,t,s,E}$,  close to the transition (see Fig.~\ref{fig_eigenvector}). The latter is very well fitted by $a(x) \propto 1/|x|$ in a wide range of $x$, except for some cutoff near $x=0$ and $x \rightarrow \pm \infty$.
The general idea of our proof is then simple: if the operator $F_{K,t,s,E}$ was acting on a real Hilbert space of finite dimension, then we could use the Collatz-Wielandt formula \cite{Meyer_2000} to derive \beq \label{eq_cw_finite_d} \sup_{v \succeq 0} \inf_{\substack{x \\ v(x) >0}} \frac{\int \mathcal{K}_{K,t,s,E}(x,y) v(y) \dd y}{v(x)}   = \lambda_{K,t,s,E} = \inf_{u \succ 0} \sup_{x} \frac{\int \mathcal{K}_{K,t,s,E}(x,y) u(y) \dd y}{u(x)} \ . \eeq
Here the order is the partial order defined by $u \preceq v$ iff $u(x) \leq v(x)$ for (almost) all $x$. Then, an asymptotic study of the conditions (\ref{eq_cw_finite_d}) for $K$ large, using well-chosen test vectors, would allow us to conclude. However, $F_{K,t,s,E}$ acts on a Banach space of infinite dimension, and is not compact (if it was, we could use natural extensions of finite dimension spectral results \cite{book_schaefer}). We shall circumvent this difficulty in the next section by using an elementary version of the Collatz-Wielandt formula. Note that the fact that $\rho_{K,t,E}$ appears in $\mathcal{K}_{K,t,s,e}$, instead of the bare density $\rho$, should not be relevant in the large $K$ limit, as was already noted in the early works \cite{ACTA73,Thouless74} -- even though it will add some technical difficulties in the following.

Finally, our method is not (in its present form) efficient for $s < 1$, but one can still study Eq. (\ref{eq_ACTA}) numerically to estimate the free energy  function $\varphi_{K,t}(s,E)$. What one finds is that $\varphi_{K,t}(s,E)$ is equivalent to $s \log t$ for $s$ going to $0$ (as was noticed in \cite{aw_prl_11} for the Cauchy disorder case). Such a linear approximation to $\varphi_{K,t}(s,E)$ would lead to a critical value of the hopping scaling like $1/K$. However, for $s>0$, there are subleading (and disorder dependent) corrections to this linear form that go as $\log (-\log t)$ for small $t$, giving rise to the correct scaling $t_K = g_c /(K \log K)$ for the critical value of the hopping strength.

\section{Large $K$ asymptotics for $\varphi_{K,t}(s,E)$ and $\varphi_{K,t}^{(\LL)}(s,E,L)$}
\label{sec_bounds}

\subsection{Elementary bounds on the rate of growth of a positive kernel}
As explained above, our proof is essentially based on the following elementary proposition:

\begin{myprop}Let $(X,\| \cdot \|)$ be a normed vector space with a partial order $\preceq$ compatible with the product by a non-negative real ($u \preceq v$ and $\lambda \geq0 \Rightarrow \lambda u \preceq \lambda v $), and such that \beq \forall u,v \in X,Ê\; (0 \preceq u \preceq v) \Rightarrow \|u\| \leq \|v\| \ . \eeq Let $F$ be a linear application from $X$ to itself preserving $\preceq$ ($u \preceq v \Rightarrow F(u) \preceq F(v)$). Let $u$ and $v$ be non-negative vectors ($0 \preceq u, 0 \preceq v$), different from the null vector, satisfying respectively: \beq \label{eq_weak_cw_hyp1} \begin{split}  F(u) &\preceq \lambda \; u \ , \\   F(v) & \succeq \mu\; v \end{split} \ , \eeq
for some $\lambda,\mu >0$.
Then for all vector $p \in X$ such that there exist $a,b>0$ for which \beq \label{eq_weak_cw_hyp2} a v \preceq p \preceq b u \ , \eeq it holds that:
\beq \label{eq_weak_cw} \mu \leq \liminf_{n \rightarrow \infty}  \| F^n(p) \|^{1/n} \leq  \limsup_{n \rightarrow \infty}  \| F^n(p) \|^{1/n} \leq \lambda \ . \eeq
\end{myprop}
The proposition easily follows from (\ref{eq_weak_cw_hyp1}) and (\ref{eq_weak_cw_hyp2}) by successive applications of the operator $F$.
\\

In the following, without loss of generality we prove Lemma~\ref{th_precise_1} and~\ref{th_precise_2} for $E=0$ -- the generic case can always be recovered by considering the shifted density $\rho(\cdot + E)$, which satisfies the same hypotheses \textbf{(A-C)} as $\rho$. We will apply Proposition 1 to the operators $F_{K,t,s,E}$ and $F^{(\LL)}_{K,t,s,E}$ on the normed vector space $L^1(\mathbb{R})$, with the partial order $\preceq$ defined by $u \preceq v$ if and only if $u(x) \leq v(x)$ for almost all $x \in \mathbb{R}$. Since $F_{K,t,s,E}$ and $F^{(\LL)}_{K,t,s,E}$ are linear and have a non-negative kernel, it is an easy check that they indeed preserve the order $\preceq$.

\subsection{Lower bound on $ \varphi_{K,t}(s,E)$ in the localized phase}
\label{sec_lower_bound}
In order to prove Lemma~\ref{th_precise_1}, we assume that $E=0$ lies in the localized phase, meaning that $(K,t) \in \textrm{Loc}$ where
\beq \begin{split} \textrm{Loc} = \{(K,t), \textrm{there exists a neighborhood $\V$ of $E=0$ such that the spectrum of $H_{K,t}$}\\\textrm{ exhibits spectral localization in $\V$.}\} \end{split} \eeq
Because of the technical reasons mentioned in Sec.~\ref{sec_comp_loc}, we still have to consider a small interval of energies around $E=0$. 
For $(K,t) \in \textrm{Loc}$ we fix $0 < \delta_{K,t} <1/(2K)$ such that the spectrum of $H_{K,t}$ is localized within $[-\delta_{K,t},\delta_{K,t}]$, and we define \beq \label{eq_def_Ecal} \mathcal{E}_{K,t} = \{ E \in [-\delta_{K,t},\delta_{K,t}], \; \lim_{\eta \searrow 0} \Gamma_{K,t}(y,E+i\eta,\omega) \textrm{ exists and is real for all site $y$ and almost all $\omega$} \} \ . \eeq
This allows us to use formula (\ref{eq_varphi_final}) to compute $\varphi_{K,t}(s,E)$ for $E \in \mathcal{E}_{K,t}$, while the bound on $\delta_{K,t}$ is needed mostly for technical reasons. As already mentioned, we know that $[-\delta_{K,t},\delta_{K,t}] \setminus \mathcal{E}_{K,t}$ has zero Lebesgue measure (and we conjecture that this is actually an empty set). Finally, we also fix $G>0$ such that, for all $K$, the hopping $t$ belongs to the interval $[0,G/K]$.
\\

We consider for $(K,t) \in \textrm{Loc}$ and $E \in \mathcal{E}_{K,t}$ fixed, the following ``test-vector'' $v$, for a given $\Delta \in (0,1)$:
\beq v(x) =  \left\{ \begin{array}{l l} \frac{1}{|x|}\hspace{1 cm}  &\textrm{if } \Delta \leq |x| \leq 1\\ 0  & \textrm{otherwise.} \end{array} \right .  \eeq
As discussed in Sec.~\ref{sec_connection}, for $\Delta$ well-chosen, this vector is very close to a numerically determined dominant eigenvector for $|x| \in (0,1)$; it will turn out that the expected shape of the eigenvector for $|x| >1$ is not needed for this variational computation.

We shall prove in the appendix (cf Eq. (\ref{eq_app_positivity})) that the vector $p_{K,t,E}(-\Gamma/t^2)$ is uniformly bounded away from zero on any interval of $\mathbb{R}^*$. Since $v(x)$ has bounded support and is bounded, (\ref{eq_weak_cw_hyp2}) will then be satisfied for some $a>0$. Moreover, $v(x)$ satisfies $F_{K,t,s,E}(v) \succeq \mu v$ as soon as:
\beq \label{eq_CN_mu} \mu \leq \inf_{|x| \in [\Delta,1]} \frac{ \int \mathcal{K}_{K,t,s,E}(x,y) v(y) \dd y}{v(x)} \ . \eeq
Hence \beq \label{eq_lb_abstract_varphi} \varphi_{K,t}(s,E) \geq \log \left[ \inf_{|x| \in [\Delta,1]} \frac{ \int \mathcal{K}_{K,t,s,E}(x,y) v(y) \dd y}{v(x)} \right] \ . \eeq
We compute, for $x \in [\Delta,1]$:
\beq \label{eq_lower_bound_kernel} \begin{split} \frac{\int \mathcal{K}_{K,t,s,E}(x,y) v(y) \dd y}{v(x)} &= \frac{t^{2-s}}{|x|^{1-s}} \int_{\Delta \leq |y| \leq 1} \rho_{K,t,E} \left(-y -\frac{t^2}{x} \right) \frac{1}{|y|} \dd y \ .
 \end{split} \eeq
The idea is now to take $\Delta$ small when $K \rightarrow \infty$, so that the above integral becomes dominated by what happens for $y$ small. Since moreover $\rho_{K,t,E}$ should resemble $\rho(\cdot+E)$ for $K$ large and $t$ accordingly small, the integral in (\ref{eq_lower_bound_kernel}) should become equivalent to $-2 \rho(E-t^2/x) \log \Delta$. This can be proved by noting that $\rho_{K,t,E}$ share the same Lipschitz continuity property as $\rho$ (this is proved in the appendix, cf. Eq. (\ref{eq_rho_lip})). We therefore obtain the following bound,
for all $K \geq 2$, $t \leq G/K$, $E \in \mathcal{E}_{K,t}$ and $z \in \mathbb{R}$: \beq \label{eq_new_lemma1} -2 \rho_{K,t,E}(z)\log \Delta - 2 C \leq \int_{\Delta \leq |y| \leq 1} \rho_{K,t,E} \left(-y +zÊ\right) \frac{1}{|y|} \dd y \leq  -2 \rho_{K,t,E}(z) \log \Delta + 2C \ . \eeq
 Replacing in (\ref{eq_lower_bound_kernel}), we obtain that:
 \beq \begin{split} \frac{\int \mathcal{K}_{K,t,s,E}(x,y) v(y) \dd y}{v(x)}  \geq \frac{t^{2-s}}{|x|^{1-s}} \left[-2 \rho_{K,t,E} \left(- \frac{t^2}{x} \right)  \log \Delta - 2 C \right] \ .
 \end{split} \eeq
 This leads, using (\ref{eq_lb_abstract_varphi}), taking $s \nearrow 1$ and setting $\Delta = t^2/\alpha$, with $\alpha \in (0,1)$, to:
  \beq \label{eq_lower_bound_last} \varphi_{K,t}(1,E) \geq \log t + \log \left[ -4 {m}_{K,t}(\alpha)  \log t + 2   {M}_{K,t}(\alpha) \log \alpha - 2 C  \right]_+ \ , \eeq
where $[x]_+ = \max(x,0)$ and we introduced
 \beq {m}_{K,t}(\alpha) = \inf_{\substack{|z|\leq \alpha\\ E \in \mathcal{E}_{K,t}}} {\rho}_{K,t,E}(z) \ , \hspace{2 cm} {M}_{K,t}(\alpha) = \sup_{\substack{|z| \leq \alpha\\ E \in \mathcal{E}_{K,t}}} {\rho}_{K,t,E}(z) \ . \eeq
Remembering that $\rho_{K,t,E}$ is the density of $V-E-t^2  \sum_{i=1}^{K-1} \Gamma_i$, in the limit of $K$ large but with $t \leq G/K$ it should get close to the bare density $\rho(\cdot +E)$. More formally, we prove in the appendix (cf Eq. (\ref{eq_cv_large_K_inter})) that, in the limit of $K$ large, with $t_K$ depending on $K$ and satisfying $t_K \leq G/K$:
\begin{align} \label{eq_limit_mk} \lim_{K \rightarrow \infty} {m}_{K,t_K}(\alpha) &= m(\alpha) \equiv   \inf_{|z| \leq \alpha}\rho(z) \ ,  \\ \label{eq_limit_Mk} \lim_{K \rightarrow \infty} {M}_{K,t_K}(\alpha) &= M(\alpha) \equiv \sup_{|z| \leq \alpha} \rho(z) \ .  \end{align}
From (\ref{eq_lower_bound_last}) and (\ref{eq_limit_mk}-\ref{eq_limit_Mk}) we deduce the asymptotic behavior of $\varphi_{K,t_K}(1,E)$ near $E=0$ when $t_K = \frac{g}{K \log K}$ and $g$ is large enough:
\beq \forall g > \frac{1}{4  m(\alpha)}, \;\;\;\; \liminf_{\substack{K \rightarrow \infty \\ \left(K,\frac{g}{K \log K} \right) \in \textrm{Loc}}}  \essinf_{E \in [-\delta_{K,t_K},\delta_{K,t_K}]} \varphi_{K,\frac{g}{K \log K}} (1,E) + \log K  > 0 \ . \eeq
Because this is valid for $\alpha$ arbitrarily small, we deduce Lemma~\ref{th_precise_1}:
\beq \forall g > \frac{1}{4 \rho(0)}, \;\;\;\; \liminf_{\substack{K \rightarrow \infty \\ \left(K,\frac{g}{K \log K} \right) \in \textrm{Loc}}} \essinf_{E \in [-\delta_{K,t_K},\delta_{K,t_K}]}  \varphi_{K,\frac{g}{K \log K}} (1,E)  + \log K > 0 \ . \eeq

\subsection{Upper bound on $ \varphi_{K,t}^{(\LL)}(s,E,L)$}
\label{sec_upper_bound}
To derive an upper bound on $\varphi_{K,t}^{(\LL)}(s,E,L)$, we consider for $\Delta  \in (0,1)$ and $A>0$ the following vector $u \in L^1(\mathbb{R})$:
\beq u(x) =  \left\{ \begin{array}{l l} \displaystyle \frac{1}{|x|^{2-s}} \hspace{1 cm}  &\textrm{if } |x| \geq 1 \\  \displaystyle \frac{1}{|x|} \hspace{1 cm}  &\textrm{if } \Delta \leq |x| \leq 1 \\ \displaystyle \frac{A}{\Delta^{\varsigma}}\frac{1}{|x|^{1-\varsigma}}  & \textrm{if } |x| < \Delta \end{array} \right. \ . \eeq
($\varsigma$ was defined by condition \textbf{(B)} on $\rho$). For technical reasons, we also assume $|E| \leq 1/4$.
Using the bounds on $p_{K,t,E}^{(\LL)}$ derived in appendix (Eq. (\ref{eq_app_p_LL})), $u$ satisfies (\ref{eq_weak_cw_hyp2}) with respect to $p_{K,t,E}^{(\LL)}(- \cdot/t^2)$ for some $b >0$.
Then we have, for $|x|< \Delta$ (and denoting as before $\Delta = t^2/\alpha$):
\beq \begin{split} \label{eq_upper_bound_A} \frac{\int \mathcal{K}^{(\LL)}_{K,t,s,E}(x,y) u(y) \dd y}{u(x)} &= \frac{t^{2-s}}{|x|^{1-s+\varsigma}} \left [ \int_{|y| \leq \Delta} {\rho_{K,t,E}^{(\LL)}}\left( -y - \frac{t^2}{x} \right) \frac{1}{|y|^{1-\varsigma}}\dd y \right. \left. + \frac{\Delta^{\varsigma}}{A}  \int_{\Delta \leq |y| \leq 1} {\rho_{K,t,E}^{(\LL)}}\left( -y - \frac{t^2}{x} \right) \frac{1}{|y|} \dd y  \right.\\& \hspace{4 cm} + \left. \frac{\Delta^{\varsigma}}{A} \int_{|y|\geq 1} \rho_{K,t,E}^{(\LL)} \left(-y - \frac{t^2}{x}\right) \frac{1}{ |y|^{2-s}} \dd y\right] \\&\leq {t^{2-s}} \frac{\Delta^{\varsigma}}{t^{2-2s+2 \varsigma}} \sup_{|z| \geq \alpha} \left [ {2}\left( \frac{1}{\varsigma} -\frac{1}{A} \log \Delta \right) \sup_{\substack{|y| \leq 1}}  \rho_{K,t,E}^{(\LL)} \left(y+z \right) \left |  z\right|^{1-s+\varsigma}  \right. \\ & \hspace{4 cm} \left.+ \frac{1}{A}\int_{|y|\geq 1} \rho_{K,t,E}^{(\LL)} \left(y + z\right) \frac{|z|^{1-s+\varsigma}}{ |y|^{2-s}} \dd y \right] \end{split} \ . \eeq
 We now use the following lemma, proved in the appendix:
 \begin{mylemma}\label{lemma_app} There exist $\mathcal{C}_1, \mathcal{C}_2 \in \mathbb{R}$ such that, for all $K \geq 2, \LL \geq 2, t \leq G/K, |E| \leq 1/4$ and $s \geq (1+\varsigma)/2$:
 \begin{align}  \sup_{z} \sup_{\substack{|y| \leq 1}}   \rho_{K,t,E}^{(\LL)}\left(y+z \right) \left | z^{1-s+\varsigma} \right|  \leq  \mathcal{C}_1 \ , \\
  \sup_{z} \int_{|y|\geq 1} \rho_{K,t,E}^{(\LL)} \left(y + z\right) \frac{|z|^{1-s+\varsigma}}{ |y|^{2-s}}  \mathrm{d} y \leq  \mathcal{C}_2 \ . \end{align} 
  \end{mylemma}
This gives by replacing in (\ref{eq_upper_bound_A}), and for $s \geq (1+\varsigma)/2$: \beq  \label{eq_41} \sup_{|x| < \Delta} \frac{\int \mathcal{K}^{(\LL)}_{K,t,s,E}(x,y) u(y) \dd y}{u(x)} \leq t^{s} \alpha^{-\varsigma} \left [ \left( \frac{2}{\varsigma} - \frac{4}{A}\log t + \frac{2}{A} \log \alpha \right) \mathcal{C}_1 +  \frac{1}{A} \mathcal{C}_2  \right] \ .\eeq
  \\
On the other hand, for $|x| \geq \Delta$, and using again the Lipschitz continuity of $\rho^{(\LL)}_{K,t,E}$ and Eq. (\ref{eq_new_lemma1}):
\beq \begin{split} \label{eq_42} \frac{\int \mathcal{K}_{K,t,s,E}^{(\LL)}(x,y) u(y) \dd y}{u(x)} &\leq \frac{t^{2-s}}{|\Delta|^{1-s}} \left [2\frac{A}{\varsigma}  \| \rho_{K,t,E}^{(\LL)} \|_\infty+   \int_{\Delta \leq |y| \leq 1} \rho_{K,t,E}^{(\LL)} \left( -y - \frac{t^2}{x} \right) \frac{1}{|y|}\dd y + 1 \right]
\\ & \leq t^s \alpha^{1-s} \left[ 2 A \frac{\|\rho\|_\infty}{\varsigma} - 2\rho_{K,t,E}^{(\LL)} \left( -\frac{t^2}{x} \right)  \log \Delta  + 2 C+1 \right] \ .
 \end{split} \eeq
 To go from the first line to the second, we used that $\| \rho_{K,t,E}^{(\LL)}\|_\infty \leq \|\rho\|_\infty$; this follows from the fact that $\rho^{(\LL)}_{K,t,E}$ is the density of a sum of random variables that contains a random variable $V$ distributed with density $\rho$, and from the bound on the convolution of two densities $f_1,f_2$: $\| f_1 \star f_2 \|_\infty \leq \|f_1\|_\infty$.  Taking $A = 2 \mathcal{C}_1 \frac{\alpha^{-\varsigma}}{\rho(0)}$, and keeping $t \leq G/K$, (\ref{eq_42}) will be larger than (\ref{eq_41}) for $K$ large enough (independently of $s \geq (1+\varsigma)/2$). Therefore, for every $\epsilon >0$, $K$ large enough, $s$ close enough to one and $L$ large enough:
  \beq \label{eq_bound_finite_volume} \varphi_{K,t}^{(\LL)}(s,E,L) \leq \log t + \log  \left[ 4 \mathcal{C}_1 \frac{\alpha^{-\varsigma}}{\rho(0)} \frac{ \| \rho \|_\infty}{\varsigma} -4 {M}^{(\LL)}_{K,t,E}(\alpha)  \log t  + 2 {m}^{(\LL)}_{K,t,E}(\alpha) +  \log \alpha + 2 C + 1 \right]_+ + \epsilon \ . \eeq
In this equation we defined:
  \beq {m}_{K,t,E}^{(\LL)}(\alpha) = \inf_{\substack{|z|\leq \alpha }} {\rho}^{(\LL)}_{K,t,E}(z) \ , \hspace{2 cm} {M}_{K,t,E}^{(\LL)}(\alpha) = \sup_{\substack{|z| \leq \alpha}} {\rho}^{(\LL)}_{K,t,E}(z) \ . \eeq
  Moreover, the bound (\ref{eq_bound_finite_volume}) is satisfied, for $\epsilon>0$ fixed, uniformly in $K$ and $L$ when $E$ is varied; this follows from a closer look at Proposition 1, and on the fact that the coefficient $b$ appearing in Eq. (\ref{eq_weak_cw_hyp2}) can be choosen independently of $E$.   
  Again, using the convergence of $\rho_{K,t,E}^{(\LL)}$ towards $\rho$ when $t \leq G/K$ (cf. Eq. (\ref{eq_cv_large_K_inter_0}) of the Appendix, and the remark below it), we have:
 \begin{align}  \lim_{K \rightarrow \infty} {m}_{K,t,E}^{(\LL)}(\alpha) &= m_E(\alpha) \equiv   \inf_{|z| \leq \alpha}\rho(z+E) \ ,  \\  \lim_{K \rightarrow \infty} {M}_{K,t,E}^{(\LL)}(\alpha) &= M_E(\alpha) \equiv \sup_{|z| \leq \alpha} \rho(z+E) \ .  \end{align}
  Moreover, the convergence is uniform for $|E| \leq 1/4$ and $\LL \geq 2$. From this, we deduce the asymptotic upper bound:
\beq \begin{split} \forall g < \frac{1}{4 M_0(\alpha)},&\ \exists \epsilon_0 >0, \ \exists \delta >0,\ K_0 \geq 0,\ \forall K \geq K_0,\ \exists s <1,\ L_0 \geq0,\ \forall L \geq L_0,\ \forall \LL \geq 2, \\& \sup_{E \in [-\delta,\delta]} \varphi_{K,\frac{g}{K \log K}}^{(\LL)} (s,E,L)  + \log K \leq -\epsilon_0 \ .  \end{split} \eeq
Because this is valid for $\alpha$ arbitrarily small, we deduce Lemma~\ref{th_precise_2}.

\section{Numerical results}
In order to demonstrate the agreement of our results with numerical studies, we performed several numerical simulations. To simplify the discussion, we will assume in this section that there exists, for the energy $E$ considered, and for every connectivity $K+1$, a unique number $g_c(E,K)$ such that the spectrum of $H_{K,\frac{g}{K \log K}}$ is localized in a neighborhood of $E$ if $g < g_c(E,K)$ and absolutely continuous in a neighborhood of $E$ otherwise. In particular $g_c(E) = \lim_{K \rightarrow \infty} g_c(E,K)$. We will focus on results obtained for two choices of the probability density $\rho$ of the disorder: the uniform distribution $\rho(v) = \frac{1}{2} \mathbf{1}( v \in [-1,1])$ -- this distribution is obviously not Lipschitz continuous on $\mathbb{R}$, but it is around $E=0$ and we believe that our result does not apply in this case only for technical reasons in the proof -- and the Cauchy distribution $\rho(v) = \frac{1}{\pi (v^2+1)}$. In fact, in both cases, as we will explain later, we have a prediction not only for the leading order of $\varphi_{K,t}(1,E)$ within the localized phase, but also for the first subleading correction. Moreover, in the Cauchy case, it is well-known that the density $\rho_{K,t,E}$ can be computed exactly, which will simplify some computations.
\\

The simplest quantity to determine numerically is the largest eigenvalue of a finite dimensional approximation of the kernel $\mathcal{K}_{K,t,s,E}$, with the additional assumption that one may substitute $\rho$ for $\rho_{K,t,E}$ in $\mathcal{K}_{K,t,s,E}$. Since this (modified) kernel is positive and irreducible, the latter can be found by iterations -- a few are enough. One then has to let the number of discretization steps go to infinity in order to ``correctly'' approximate the original continuous operator (this remains however uncontrolled). Since we do not use the exact density of states $\rho_{K,t,E}$ in this case, this is still an approximation of the operator $F_{K,t,s,E}$ which is expected to give good results for the mobility edge $g_c(K,E)$ only as $K$ becomes large. More precisely, one expects the finite $K$ corrections due to the use of $\rho$ instead of $\rho_{K,t,E}$ to be quite small (of order $K t^2 =O(1/(K \log^2 K)$) \cite{ACTA73,Thouless74}); in particular much smaller than those in powers of $1/\log K$ that we expect  to find regardless of the precise form of $\rho$. We denote $g_c^{\mathbf{(A)}}(K,E)$ the critical value computed with this procedure.

In order to deal with the true kernel $\mathcal{K}_{K,t,s,E}$, one has to determine the modified density of states $\rho_{K,t,E}$ first. For the Cauchy distribution this can be done analytically \cite{ACTA73,MiDe94}; in general one first has to solve the recursive equation (\ref{eq_rde_gamma_real}) with some population dynamics approximation \cite{ACTA73,MiDe94,MoGa09,BST10}  -- the precision of the latter resolution is not critical in this case, and moderate population sizes and number of iterations are enough to get a good approximation of $\rho_{K,t,E}$. Plugging this into $\mathcal{K}_{K,t,s,E}$ and following the same route as before to find the largest eigenvalue of this kernel this gives a critical value $g_c^{\mathbf{(B)}}(K,E)$ which is exact (assuming infinite numerical precision) and is the most convenient to compute. It can be seen in Table~\ref{table_uniform}, \ref{table_cauchy} that this value is, as anticipated in the previous paragraph, very well approximated by $g_c^{\mathbf{(A)}}(K,E)$ even for moderate values of $K$.
\\

A different way to proceed is to look for the localization transition directly on the tree. In this case it is actually much more convenient numerically to consider the ``quenched free energy'' \beq \phi_{K,t}(s,E) \equiv \lim_{\eta \searrow 0} \lim_{R \rightarrow \infty} \frac{1}{R} \mathbb{E}  \log \left[ \sum_{\substack{x \in \mathcal{T}_K \\ \textrm{dist}(0,x) = R}} \left|G_{K,t}(0,x,E+ i\eta, \omega)\right|^s\right] - \log K \ , \eeq
(note that compared to (\ref{eq_def_varphi}), the log is outside of the sum). The latter can be computed, for any $s$, with the cavity method \cite{BTT}. The localization transition was usually found \cite{ACTA73,MiDe94,MoGa09,BST10} by looking at the critical value of $t$ such that $ \phi_{K,t}(2,E)$ equals $-\log K$, because this corresponds to the apparition of a non-vanishing imaginary part in $\Gamma$ solution of (\ref{eq_rde_gamma_im}) when $\eta \searrow 0$ \cite{aw_long_v3,aw_prl_11,BTT}. In this case, strong finite-size effects were mentioned  \cite{ACTA73,MiDe94,MoGa09,BST10,BTT} and the values obtained were hard to reconcile with those obtained with the method presented above \cite{ACTA73}. However, at the localization transition, $\phi_{K,t}(s,E)$ does not depend on $s$ for $s \in [1,2]$ and is equal to $\varphi_{K,t}(1,E)$ \cite{BTT}. One can therefore compute $\phi_{K,t}(s,E)$ as a function of $t$  for any $s \in [1,2]$ and check when it reaches the critical value $-\log K$. It turns out that doing this simulation for $s=1$ strongly reduces the finite-size effects. Moreover it also allows to take $\eta = 0$ from the beginning (within the localized phase), reducing the number of parameters on which the numerical results depend. In addition to that, we performed two finite-size scaling analysis. To explain them, let us call $\phi_{K,t}^{(R,N)}(s,E)$ the approximated value of $\phi_{K,t}(s,E)$ computed with $R$ iterations of a pool of size $N$ (the procedure is analog to that of \cite{BTT}). We first take the $R \rightarrow \infty$ limit at fixed $N$, using a finite-size correction as $\phi_{K,t}^{(R,N)}(s,E) = \phi_{K,t}^{(\infty,N)}(s,E) + a_{K,t}(N)/R$ and $R$ in the range $[1000,10000]$. Then we take the population size $N$ to infinity, assuming logarithmic corrections of the form $\phi_{K,t}^{(\infty,N)}(s,E) = \phi_{K,t}^{(\infty,\infty)}(s,E) + b_{K,t}/\log N + c_{K,t}/ (\log N)^2$. We used this form because it fits relatively well the corrections, and also because it was predicted in \cite{MoGa09} for this particular model; of course such slow corrections impoverish the quality of the numerical fit. For the purpose of this work, we used population sizes between $10^4$ and $4.10^6$, and also took an extra average over a few independent realizations of the $R$ cavity iterations, in order to reduce fluctuations effects. All of this leads to a computationally heavy procedure. However it gives critical values $g_c^{\mathbf{(C)}}(K,E)$ that are quite different from the ones previously computed with this method \cite{ACTA73,BST10}, but very close to the exact value $g_c^{\mathbf{(B)}}(K,E)$ computed with the procedure above. In particular this resolves the discrepancy that was observed in \cite{ACTA73}.
\\

Finally we want to compare these results with analytic asymptotic predictions for $g_c(K,E)$. The first one is a straightforward consequence of the bounds of Sec.~\ref{sec_bounds}, and defines $g_c^{\mathbf{(D)}}(K,E)$ as the smallest root of \beq  \log g -\log  \log K +  \log \left[ -4 \rho(E) \left( \log g - \log(K\log K) \right)\right]\ . \eeq This expansion can in principle be improved by computing the next subleading term in $\varphi_{K,t}(1,E)$ for $K$ large. In the particular case of $E=0$, one should recover the prediction of \cite{ACTA73} (Eq. (7.6)), obtained through an approximate solution of the eigenvalue equation (\ref{eq_ACTA}) for the operator $F_{K,t,s,E}$. This gives, both in the uniform and in the Cauchy case, a critical value $g_c^{\mathbf{(E)}}(K,0)$ defined as the smallest root of: \beq \log g - \log \log K + \log \left[- 4 \rho(0) \left( 1+ \frac{\pi^2}{24} \frac{1}{(\log g - \log (K \log K) )^2} \right)\left( \log g - \log(K\log K) \right) \right] \ . \eeq Note that there is no reason to expect this expansion, that may seem universal at first sight (the same numerical factor $\pi^2/24$ appearing in both cases), to remain valid for other disorder densities, or for $E \neq 0$. Also note that in the uniform case, this last result can be recovered through a systematic expansion of $\varphi_{K,t}(s,E)$ in powers of $1/\log K$, starting from Eq. (\ref{eq_G_unfolded_final}) \cite{BaMu_unp}. In any case, we expect that \cite{footnote_2}: \beq g_c^{\mathbf{(D)}}(K,0) = g_c(K,0) + O \left( \frac{1}{(\log K)^2}\right) \ , \hspace{1cm} g_c^{\mathbf{(E)}}(K,0) = g_c(K,0) + O \left( \frac{1}{(\log K)^3}\right) \ , \eeq
where the $O(\dots)$ terms should be rather small, even at finite $K$. On the other hand we expect for the critical value at finite $K$, $g_c(K,E)$, a very slow convergence towards its limiting value: \beq g_c(K,E) = g_c(E)- \frac{1}{4 \rho(E)} \frac{ \log \log K}{\log K} + O \left( \frac{1}{\log K} \right)  \ . \eeq
Hence it is worth considering the subleading corrections that led to $g_c^{\mathbf{(D,E)}}(K,E)$.
The agreement between these values (in particular $g_c^{\mathbf{(E)}}(K,0)$) and the best numerical value $g_c^{\mathbf{(B)}}(K,0)$ can be seen to be very good in the uniform case, and also satisfying in the Cauchy case, even though the finite $K$ corrections are stronger.

 \label{sec_numerical}
 \begin{table}[h]
\begin{tabular}{|c||c|c|c||c|c|}
\hline
$K$ & $g_c^{\mathbf{(A)}}(K,0)  $ & $g_c^{\mathbf{(B)}}(K,0) $ & $g_c^{\mathbf{(C)}}(K,0) $ &  $g_c^{\mathbf{(D)}}(K,0) $ &  $g_c^{\mathbf{(E)}}(K,0) $\\ \hline
 2 & 0.150 & 0.153 &0.154 &0.154 &0.149 \\
 3 & 0.187 & 0.188 &0.189 &0.194 &0.187 \\
 4 & 0.207 & 0.208 &0.204 &0.213 & 0.207\\
 5 & 0.220 & 0.220 &0.219 &0.225 & 0.220\\
 6 & 0.230 & 0.231 & 0.227 &0.234 & 0.230 \\
 8 & 0.243 & 0.243 &- &0.247 & 0.243\\
 12 & 0.261 & 0.261 &- & 0.263 &  0.260\\
 \hline
\end{tabular}
\caption{Critical values of the disorder for the uniform case and $E=0$, for several numerical procedures (superscript \textbf{A} to \textbf{C}) and two asymptotic formula (superscript \textbf{D} and \textbf{E}), defined in the text. The numerical error for $g_c^{\mathbf{(A)}}$ and $g_c^{\mathbf{(B)}}$ is smaller than $0.001$, while it is larger for $g_c^{\mathbf{(C)}}$ (of the order of $0.005$). In this case the asymptotic value is $g_c(E)= 1/2$, and is approached extremely slowly. The connection with the model where $V$ is uniformly distributed in $[-W/2,W/2]$ and $t=1$ is given by $W_c=\frac{2}{g_c} K \log K$.}
\label{table_uniform}
\end{table}

 \begin{table}[h]
\begin{tabular}{|c||c|c|c||c|c|}
\hline
$K$ & $g_c^{\mathbf{(A)}}(K,0)  $ & $g_c^{\mathbf{(B)}}(K,0) $ & $g_c^{\mathbf{(C)}}(K,0) $ &  $g_c^{\mathbf{(D)}}(K,0) $ &  $g_c^{\mathbf{(E)}}(K,0) $
\\ \hline
 2 & 0.317 & 0.334 &0.334 &- &0.367 \\
 3 & 0.364 & 0.372 &0.370 &0.418 &0.384 \\
 4 & 0.389 & 0.394 &0.394 &0.423 &0.403 \\
 5 & 0.406 & 0.410 &0.404 &0.432 &0.417 \\
 6 & 0.419 & 0.421 & 0.422 &0.440 & 0.428 \\
 8 & 0.436 & 0.437 &- &0.453 &0.444  \\
 12 & 0.456 & 0.457 &- & 0.470& 0.463 \\
 \hline 
\end{tabular}
\caption{Critical values of the disorder for the Cauchy case and $E=0$, for several numerical procedures (superscript \textbf{A} to \textbf{C}) and two asymptotic formula (superscript \textbf{D} and \textbf{E}), defined in the text. The numerical error for $g_c^{\mathbf{(A)}}$ and $g_c^{\mathbf{(B)}}$ is smaller than $0.001$, while it is of order $0.01$ for $g_c^{\mathbf{(C)}}$. In this case the asymptotic value is $g_c(E)= \pi/4 \simeq 0.785$. The connection with the model where $V$ has Cauchy distribution of parameter $\gamma$ and $t=1$ is given by $\gamma_c=\frac{1}{g_c} K \log K$.}
\label{table_cauchy}
\end{table}

\section{Conclusion}
\label{sec_ccl}
In this article we have computed and bounded the free energy function $\varphi_{K,t}$ introduced in \cite{aw_long_v3} on tree graphs and in the large connectivity limit, either within the pure point phase or in a finite volume setting. These two bounds allow to elucidate the asymptotic scaling of the mobility edge in this large connectivity regime.
Interestingly, we have found that this rigorous approach gives back the criterion obtained in \cite{ACTA73}
using a self consistent equation on the tail of the distribution of the local Green function.

Our technique applies at any energy $E$ at which the original density of disorder $\rho(E)$ is continuous
and strictly positive (and under additional assumptions on $\rho$ that are probably mostly technical).
An interesting extension of this work would be to understand what happens when $\rho(E) = 0$. This
can happen in two ways: a possibility is that the energy $E$ considered is such that $\rho(E) = 0$ but
some (one-sided) derivative of $\rho$ does not vanish. Another appealing regime is when $E$ is outside
the support of the disorder density $\rho$. In this case, one may expect, following the bound given in
\cite{AW11}, a mobility edge scaling like $1/K$ in a range of energies at distance at most $1/\sqrt{K}$ from the
support of $\rho$. The behavior of the free energy function $\varphi_{K,t}$ and of the related localization length
in this Lifshitz tail regime would be particularly interesting to understand.

Besides the value of the mobility edge, another direction of investigation concerns the nature of
the delocalization transition on the Bethe lattice. It would be nice to understand the behavior of the free energy function $\varphi_{K,t}$ in the delocalized phase but close to the transition; this could in particular allow to prove the presence of only absolutely continuous spectrum. This would probably require to analyze the appearance of a small imaginary part in the local resolvent near the transition, an analysis that was non-rigorously carried in \cite{MiFy91}.

Finally, an extension of the lower bound on the localization threshold to the Anderson model on the hypercube, or to the finite dimensional Anderson model, are also stimulating points. In this last case, it is to be expected that the best known rigorous bound \cite{Schenker13} is again off by a factor $2/e$ (see \cite{Efetov88} for a non-rigorous discussion of a related model).

\acknowledgments
I am deeply indebted to M. M\"{u}ller and G. Semerjian for numerous, helpful and stimulating discussions related to this work, and for their constant support. I am particularly grateful to G. Semerjian for his careful reading of the manuscript, and to S. Warzel for inspiring discussions on a former version of it. I also thank G. Biroli, C. Bordenave, M. Lelarge and J. Salez fur useful comments on this work.

\appendix
\section{Study of the real RDE}

This appendix is devoted to the derivation of various estimates for the probability density $\rho_{K,t,E}$ in the regime of $K$ large and $t$ accordingly small ($t \leq G/K$). Our starting point is the real recursive distributional equation (\ref{eq_rde_gamma_real}), that we recall here for convenience:
\beq \label{eq_rde_real_app} \Gamma \=d \frac{1}{V-E-t^2 \sum_{i=1}^{K} \Gamma_i} \ . \eeq
In the following we shall always assume that this equation admits at least one solution, and speak about ``the solution'' of  (\ref{eq_rde_real_app}) for any given solution of it. This is not an issue since in the core of the text we use this equation only for energies $E \in \mathcal{E}_{K,t}$ (defined in Sec.~\ref{sec_bounds}), which guarantees existence of at least one solution to (\ref{eq_rde_real_app}), while the choice of the relevant solution is made unambiguous by the unicity of the solution to the complex equation (\ref{eq_rde_gamma_im}), as explained in Sec.~\ref{sec_comp_loc}. Finally, recall that we always assume in this case $E \in [-\delta_{K,t}, \delta_{K,t}]$ with $\delta_{K,t} \leq 1/(2K)$.

In order to study the finite-volume function $\varphi_{K,t}^{(\LL)}$, we need to study the random variable $\Gamma^{(\LL)}$ obtained by a finite number $\LL$ of iterations of Eq. (\ref{eq_rde_real_app}), starting with $\Gamma^{(1)} \=d \frac{1}{V-E}$, as explained in Sec.~\ref{sec_finite_volume}. In this second case, the energy $E$ is assumed to satisfy the weaker condition $E \in [-\delta_{K,t}, \delta_{K,t}]$ with $\delta_{K,t} =1/4$.
\\

We divide the results of this appendix in three  parts: first the study of the solution of equation (\ref{eq_rde_real_app}), then its consequences on the regularity of $\rho^{(\LL)}_{K,t,E}$ and the proof of Lemma~\ref{lemma_app}, and finally a study of the convergence of $\rho_{K,t,E}$ and $\rho^{(\LL)}_{K,t,E}$ towards $\rho$ when $K$ is large (and $t$ small).

\subsubsection{Probability distribution of the solution of the real recursive distributional equation}
 In this first part we derive estimates on $p_{K,t,E}$ and $p_{K,t,E}^{(\LL)}$ that were needed to apply Proposition 1 in Sec.~\ref{sec_bounds}. Here we first focus on a solution $\Gamma$ of the full recursive distributional equation (\ref{eq_rde_real_app}), before stating the straightforward generalization to the distribution of $\Gamma^{(\LL)}$.

We recall that we denoted $p_{K,t,E}$ the probability density of $\Gamma$ solution of equation (\ref{eq_rde_real_app}). Let us also denote $p^{(K)}_{K,t,E}$ (resp. $p^{(K-1)}_{K,t,E}$) the probability density of $t^2 \sum_{i=1}^{K} \Gamma_i $ (resp. $t^2 \sum_{i=1}^{K-1} \Gamma_i$). Finally, recall that we introduced the density $\rho_{K,t,E}$ of $V-E-t^2 \sum_{i=1}^{K-1} \Gamma_i $.
\\

Since $V$ admits a density, any sum of $V$ with another random variable also does \cite{feller1}. In particular this holds for the denominator in (\ref{eq_rde_real_app}), and thus also for the left hand side of (\ref{eq_rde_real_app}). Hence $p_{K,t,E}$ is indeed a density. Moreover it satisfies: \beq \label{eq_app_upper_bound_gamma_1} \begin{split} p_{K,t,E}(\Gamma) &= \frac{1}{\Gamma^2} \int p^{(K)}_{K,t,E} \left( x \right) \rho \left( \frac{1}{\Gamma} +E+ x\right) \dd x \leq \frac{\| \rho \|_\infty}{\Gamma^2} \ . \end{split} \eeq
Using that, given $n$ densities $f_1, f_2, \dots, f_n$  their convolution satisfies
\begin{align} \label{eq_bound_convolution} 
(f_1 \star f_2 \star \dots \star f_n)(x) \leq \sup_{y \geq x/n} f_1(y) + \dots +  \sup_{y \geq x/n} f_n(y)  \end{align}
(and the same inequality with $y \geq x/n$ replaced by $y \leq x/n$), we deduce that
\begin{align} \label{eq_bound_pK} p^{(K)}_{K,t,E} \left (x\right) &\leq  K^3 t^2 \frac{\| \rho \|_\infty}{|x|^2} \ . \end{align}
From assumption \textbf{(B)} on $\rho$: $\forall v \neq 0, \rho(v) \leq C_\varsigma /|v|^{1+\varsigma}$, we deduce the following bound for $x \geq 1/2$:
\beq \label{eq_useful_bound} \sup_{E \in [-\delta_{K,t},\delta_{K,t}]} \rho(x+E) \leq \sup_{E \in [-1/4,1/4]} \frac{C_\varsigma}{(x+E)^{1+\varsigma}} \leq \frac{2^{1+\varsigma}C_\varsigma}{x^{1+\varsigma}} \ . \eeq
Introducing  the density $\rhot_{K,t,E}$ of $V-E-t^2 \sum_{i=1}^{K} \Gamma_i $, and using the above equation with (\ref{eq_bound_convolution}) gives that, for $x \geq 1$
\beq  \widetilde{\rho}_{K,t,E} \left (x\right)  \leq 4 \left(2^{1 + \varsigma} \Csigma+K^3 t^2 \|\rho\|_\infty \right) \frac{1}{|x|^{1+\varsigma}}  \ . \eeq
Going back to $\Gamma$ we obtain that for $|\Gamma| \leq 1$:
\beq \label{eq_app_upper_bound_gamma_2} p_{K,t,E}(\Gamma) = \frac{1}{\Gamma^2} \widetilde{\rho}_{K,t,E} \left(\frac{1}{\Gamma}\right) \leq 4 \left(2^{1+\varsigma} \Csigma+K^3t^2 \| \rho \|_\infty \right) \frac{1}{|\Gamma|^{1-\varsigma}} \ . \eeq

The derivation of Eq. (\ref{eq_app_upper_bound_gamma_1}) only requires one iteration of Eq. (\ref{eq_rde_real_app}), and the one of Eq. (\ref{eq_app_upper_bound_gamma_2}) two iterations. This shows that their content can be extended to a result on $p_{K,t,E}^{(\LL)}$ as soon as $\LL \geq 2$, namely:
\beq \label{eq_app_p_LL} p_{K,t,E}^{(\LL)}(\Gamma) \leq \frac{\|\rho\|_\infty}{\Gamma} \ , \hspace{0.5 cm} p_{K,t,E}^{(\LL)}(\Gamma)  \leq 4 \left(2^{1+\varsigma} \Csigma+K^3t^2 \| \rho \|_\infty \right) \frac{1}{|\Gamma|^{1-\varsigma}}  \ . \eeq

We also need the positivity of $p_{K,t,E}(\Gamma)$ to derive the lower bound of Sec.~\ref{sec_lower_bound}. To do this, first note that $\Gamma \rightarrow p_{K,t,E}(\Gamma)$ is a continuous function (except possibly in zero); this follows from the Lipschitz continuity of $\rhot_{K,t,E}$, that can be proved exactly as the one of $\rho_{K,t,E}$ (cf Eq. (\ref{eq_rho_lip}) below). Hence there exists $[\alpha_-,\alpha_+] \subset \mathbb{R}$ such that the restriction of $p_{K,t,E}$ to $[\alpha_-, \alpha_+]$ is positive. Assume first that $0 \in [\alpha_-,\alpha_+]$. Since there also exists $\delta > 0$ such that the restriction of $\rho(\cdot+E)$ to $[-\delta,\delta]$ is positive, it follows from one iteration of (\ref{eq_rde_real_app}) that $\Gamma$ has a positive density near $\pm \infty$: there exists $A>0$ such that $p_{K,t,E}$ is strictly positive on $(-\infty, A) \cup (A,\infty)$. Hence $\sum_i \Gamma_i$ has a positive density on $\mathbb{R}$ (every real number can be written as the sum of $K-1$ arbitrarily large real numbers), and with another application of (\ref{eq_rde_real_app}), we deduce that the same holds for $\Gamma$. Using the uniform continuity of $p_{K,t,E}$ on every compact of $\mathbb{R}^*$, we deduce that
\beq \label{eq_app_positivity} \forall a,b \in \mathbb{R}, ab>0 \Rightarrow \inf_{\Gamma \in [a,b]} p_{K,t,E}(\Gamma) >0 \ . \eeq
It remains to prove that one can assume $0 \in [\alpha_-,\alpha_+]$. To do that, consider the successive images
of  $\alpha >0$ by the map $f : x \rightarrow \frac{1}{\epsilon-t^2 K x}$. For $\epsilon \in [-\delta, \delta]$, if $\alpha \in \textrm{supp}(p_{K,t,E})$, $f(\alpha) \in \textrm{supp}(p_{K,t,E})$. Moreover, for $\epsilon$ small enough, there exists $n$ such that $f^{2n}(\alpha) \in [-\delta/(t^2K), \delta/(t^2K)]$. Hence if $\alpha_- >0$, it follows from repeated applications of (\ref{eq_rde_real_app}) that $\sum_i \Gamma_i$ has a positive density in some subset of $[-\delta,\delta]$. The same reasoning can be done if $\alpha_+<0$, taking $\epsilon <0$ close enough to zero. Hence $\Gamma$ has a positive density near $\pm \infty$, and applying once again the recursive equation, we deduce that $p_{K,t,E}$ has a positive density near $0$, which ends the proof.
\\

We finally prove that $\rho_{K,t,E}$ is Lipschitz continuous, with the same Lipschitz constant as $\rho$. Indeed, by a straightforward computation: 
\beq \begin{split} \label{eq_rho_lip} \left | \rho_{K,t,E}(e)-\rho_{K,t,E}(e') \right| &\leq \int \left| \rho(e+x+E)-\rho(e'+x+E)\right| p^{(K)}_{K,t,E}(x) \dd x 
\\ &\leq C |e-e'| \int p^{(K)}_{K,t,E}(x) \dd x
\\ & \leq C |e-e'| \ .
\end{split} \eeq
Again, the same immediately holds for $\rho_{K,t,E}^{(\LL)}$.

\subsubsection{Proof of Lemma~\ref{lemma_app}}
We now turn to the proof of Lemma~\ref{lemma_app}. Recall that from now on it is assumed that $t \leq G/K$ with $G >0$ fixed. Here we concentrate on $p_{K,t,E}^{(\LL)}$.
\\

It will be convenient to use the following union bound, for $a<0$:
\beq \begin{split} \mathbb{P} \left( t^2 \sum_{i=1}^{K-1} \Gamma_i^{(\LL)}  < a \right) & \leq \mathbb{P} \left[ \left(\Gamma_1^{(\LL)} \leq \frac{a}{t^2(K-1)} \right) \cup \dots \cup  \left(\Gamma_{K-1}^{(\LL)} \leq \frac{a}{t^2(K-1)} \right) \right] \\ &\leq (K-1) \mathbb{P} \left( \Gamma^{(\LL)} < \frac{a}{t^2(K-1)} \right) \\ & \leq \frac{G^2 \| \rho\|_\infty}{|a|}\ , \end{split} \eeq
where we used, from (\ref{eq_app_upper_bound_gamma_1}), and for $b < 0$:
\beq \label{eq_bound_P} \mathbb{P}\left(\Gamma^{(\LL)} < \frac{b}{t^2} \right) \leq t^2 \frac{\|\rho\|_\infty}{ |b|} \leq \frac{G^2}{K^2} \frac{\|\rho\|_\infty}{ |b|} \ . \eeq
From this, and with the help of (\ref{eq_useful_bound}), we deduce that for $z > 1$
\beq \begin{split} \rho_{K,t,E}^{(\LL)}(z) = \int \rho(z+z'+E) p^{(K-1),(\LL)}_{K,t,E}(z') \dd z' &\leq \sup_{\substack{E \in [-\delta_{K,t},\delta_{K,t}]\\ z' \geq -z/2}} \rho(z+z'+E)+ \| \rho\|_\infty \mathbb{P}\left(t^2 \sum_{i=1}^{K-1} \Gamma_i^{(\LL)} < -z/2 \right) \\& \leq  \frac{4^{1+\varsigma} C_\varsigma}{z^{1+\varsigma}} + \frac{2 G^2 \| \rho\|_\infty^2}{z}\ . \end{split}  \eeq
Hence \beq  \label{eq_A12} \sup_{z \geq 2} \sup_{\substack{|y| \leq 1}}   \rho_{K,t,E}^{(\LL)}\left(y+z \right)  |z|^{1-s+\varsigma} \leq \sup_{z \geq 1} \rho_{K,t,E}^{(\LL)}(z) (z+1)^{1-s+\varsigma} \leq 8^{1+\varsigma} C_\varsigma + 4 G^2 \|\rho\|_\infty^2 \ . \eeq
The same bound holds for $z \leq -2$, while for $z \in [-2,2]$ it is clear that $ \sup_{|z| \leq 2} \sup_{\substack{|y| \leq 1}} \rho_{K,t,E}^{(\LL)}\left(y+z \right)  z^{1-s+\varsigma} \leq \|\rho\|_\infty 2^{1+\varsigma}$ (note that $\|\rho_{K,t,E}^{(\LL)}\|_\infty \leq \|\rho\|_\infty$; as explained in Sec.~\ref{sec_upper_bound}, this follows from the fact that $\rho_{K,t,E}^{(\LL)}$ is the density of a sum of random variables that contains a random variable $V$ distributed with density $\rho$, and from the bound on the convolution of two densities $f_1,f_2$: $\| f_1 \star f_2 \|_\infty \leq \| f_1 \|_\infty$). Hence we have proved the first part of the lemma.
\\

For the second part we proceed in the same way: considering first $|z| \geq 2$, we write:
\beq \begin{split}  \int_{|y| \geq 1} \rho_{K,t,E}^{(\LL)} (y+z) \frac{|z|^{1-s+\varsigma}}{|y|^{2-s}}\dd y& = \int_{\substack{|y| \geq 1\\ |y| \geq z/2}} \rho_{K,t,E}^{(\LL)} (y+z) \frac{|z|^{1-s+\varsigma}}{|y|^{2-s}}\dd y+\int_{\substack{|y| \geq 1\\ |y| \leq z/2}} \rho_{K,t,E}^{(\LL)} (y+z) \frac{|z|^{1-s+\varsigma}}{|y|^{2-s}}\dd y 
\\ & \leq 2^{2-s} |z|^{\varsigma-1} + \int_{\substack{ |y| \geq 1 \\ y \leq |z|/2}}  \left( 4^{1+\varsigma}\frac{C_\varsigma}{|z/2|^{1+\varsigma}} + 2 \frac{G^2 \| \rho\|_\infty^2}{|z/2|} \right) \frac{|z|^{1-s+\varsigma}}{|y|}\dd y
\\ & \leq 2^{2-s} |z|^{\varsigma-1} + 8 \left( 2^{1+3\varsigma} C_\varsigma +  G^2 \| \rho\|_\infty^2 \right) \frac{1}{|z|^{(1+\varsigma)/2}}  \log \left( |z|/2 \right) \ ,
\end{split} \eeq
where we used in the last line that Lemma~\ref{lemma_app} assumes $ s \geq (1+\varsigma)/2$. Hence $z \rightarrow \int_{|y| \geq 1} \rho_{K,t,E}^{(\LL)} (y+z) \frac{|z|^{1-s+\varsigma}}{|y|^{2-s}}\dd y$ is bounded uniformly in $s, K$ and $t$ for $|z| \geq 2$. Since it is also clearly the case for $|z| \leq 2$, we have proved the second part of the lemma.
\\

\subsubsection{Convergence of $\rho_{K,t,E}$ and $\rho_{K,t,E}^{(\LL)}$ towards $\rho$ for $K$ large}
\label{eq_app_cv_large_K}

Finally, we shall be interested in estimates to asses that $\rho_{K,t,E}$ is ``close'' to $\rho(\cdot+E)$. The latter are based on the intuitive idea that $t^2 \sum_{i=1}^{K-1}{\Gamma_i}$ should be of order $Kt^2 \leq G^2/K$ when $K$ is large and $t \leq G/K$. This can be proved using Eq. (\ref{eq_app_upper_bound_gamma_1}): for any $K$, $\Gamma$ is stochastically dominated by a Pareto random variable $Y$ distributed according to: \beq \rho_Y(y) = \left\{ \begin{array}{l l} \frac{\|\rho\|_\infty}{y^2}\hspace{1 cm}  &\textrm{if } y \geq \|\rho\|_\infty \ , \\ 0  & \textrm{otherwise .} \end{array} \right.   \eeq
Such a random variable satisfies the following law of large numbers \cite{breiman_book}: 
\beq  \overline{Y}_n \equiv \frac{ Y_1 + \dots + Y_n}{n} {\overset{\textnormal d}{\rightarrow}} \overline{Y} \ , \eeq
where $\overline{Y}$ is distributed according to a stable law with exponent $\alpha=1$ and tail amplitude $\|\rho\|_\infty$. In particular, for all $\epsilon>0$, $a$ large enough and $n$ large enough \beq \mathbb{P}(\overline{Y}_n\geq a) \leq (1+\epsilon) \frac{ \| \rho\|_\infty}{a} \ . \eeq
Henceforth, for $t \leq G/K$ \beq \label{eq_proba_large_K} \mathbb{P} \left( t^2 \sum_{i=1}^{K-1} \Gamma_i \geq a \right) \leq \mathbb{P} \left( Y_{K-1} \geq \frac{K^2a}{G^2(K-1)}  \right) \leq \frac{1}{ab} \ , \eeq
for any $a,b>0$, and $K$ large enough (independent of $E$). Similarly for $a, b>0$ and $K$ large enough:
\beq \label{eq_proba_large_K_2} \mathbb{P} \left( t^2 \sum_{i=1}^{K-1} \Gamma_i \leq-  a \right) \leq \frac{1}{ab} \ , \eeq
This implies that $\Gamma {\overset{\textnormal d}{\rightarrow}} \frac{1}{V-E}$. More generally, one can expand the distribution of $\Gamma$ in powers of $Kt^2$, which correspond to truncate the recursive equation (\ref{eq_gamma_recursive}) after a finite number of iterations. For the purpose of our work, it will be enough to deduce from (\ref{eq_proba_large_K}-\ref{eq_proba_large_K_2}) that for all $\epsilon>0$, $a>0$ and $K$ large enough:
\begin{align} \label{eq_cv_large_K_inter_0} \forall t \leq G/K,\;\forall e,E \in \mathbb{R},\;  \inf_{e' \in [-a,a]} \rho(e+e'+E)-\epsilon \leq \rho_{K,t,E}(e) \leq   \sup_{e' \in [-a,a]} \rho(e+e'+E)+\epsilon \; . \end{align}
Indeed, one can write for $a, b>0$ and $K$ large enough:
\begin{equation} \begin{split} \rho_{K,t,E}(e) &= \int \rho(e+x+E) p^{(K-1)}_{K,t,E}(x) \textrm{d}x 
\\ & \leq \int_{-a}^a \rho(e+x+E) p^{(K-1)}_{K,t,E}(x) \textrm{d}x + \| \rho\|_\infty \int_{\mathbb{R} \setminus [-a,a]} p^{(K-1)}_{K,t,E}(x) \textrm{d}x
\\ & \leq \sup_{e' \in [-a,a]} \rho(e+e'+E) \int_{-a}^a  p^{(K-1)}_{K,t,E}(x) \textrm{d}x +2 \| \rho\|_\infty \frac{1}{ab} 
\\ & \leq \sup_{e' \in [-a,a]} \rho(e+e'+E)  +2 \| \rho\|_\infty \frac{1}{ab} \ ,  \end{split} \end{equation}
Similarly:
\begin{equation} \begin{split} \rho_{K,t,E}(e) &= \int \rho(e+x+E) p^{(K-1)}_{K,t,E}(x) \textrm{d}x 
\\ & \geq \int_{-a}^a \rho(e+x+E) p^{(K-1)}_{K,t,E}(x) \textrm{d}x 
\\ & \geq \inf_{e' \in [-a,a]} \rho(e+e'+E) \int_{-a}^a  p^{(K-1)}_{K,t,E}(x) \textrm{d}x 
\\ & \geq \inf_{e' \in [-a,a]} \rho(e+e'+E)  \left(1 - \frac{2}{ab}\right)
\\ & \geq \inf_{e' \in [-a,a]} \rho(e+e'+E)  - 2\|\rho\|_\infty \frac{1}{ab} \ .  \end{split} \end{equation}
Taking $b = \frac{\|\rho\|_\infty}{a \epsilon}$ ends the proof of (\ref{eq_cv_large_K_inter_0}). Again, the same bounds as (\ref{eq_cv_large_K_inter_0}) immediately holds with $\rho_{K,t,E}^{(\LL)}$ ($\LL \geq 2$) instead of $\rho_{K,t,E}$. 

The Lipschitz continuity of $\rho$ implies that, for $a>0$, $K$ large enough and $t \leq G/K$:
\begin{equation} \begin{split} \sup_{|E| \leq 1/K} |Ê\rho_{K, t, E}(e)-\rho(e)| &\leq \sup_{|E| \leq 1/K, e' \in [-a,a]} | \rho(e+e'+E) -\rho(e)| + \epsilon
\\ &  \leq C (a+1/K) + \epsilon \end{split} \end{equation}
In particular for $t_K \leq G/K$:
 \beq \label{eq_cv_large_K_inter} \lim_{K \rightarrow \infty} \sup_{ |E| \leq 1/K} |\rho_{K,t_K,E}(e)-\rho(e)| = 0 \ , \eeq
and the convergence is uniform with $e$.

\bibliographystyle{h-physrev}

\begin{thebibliography}{10}
\bibitem{footnote_1}
Here it would be tempting to say that this property holds for \textit{all} $E$ in the localized phase, but we do not have a proof of this statement. One can however show that any $E$ in the localized phase is almost surely not an eigenvalue of $H_{K,t}(\omega)$; this follows from the simplicity of the spectrum within the pure point phase \cite{simon94} and ergodicity (\cite{carmona}, in particular Prop.~V.2.8).
\bibitem{footnote_2}
In this particular case only; for general densities $\rho$ and energies $E$, we expect the first correction to $g_c^{\mathbf{(D)}}$ to go as $1/\log K$.
\end{thebibliography}

\end{document}